\documentclass{segabs}

\usepackage{amsmath}

\begin{document}

\title{Verification of reciprocity in anisotropic poroelastic wave simulation using symmetric Strang splitting}
\author{Morten Jakobsen$^{1}$ and Jos\'e M.\ Carcione$^{2}$\\
$^{1}$VISTA CSD, Department of Earth Science, University of Bergen \quad
$^{2}$OGS Trieste
}

\footer{Jakobsen \& Carcione}
\lefthead{Jakobsen \& Carcione}
\righthead{Reciprocity‑Preserving Biot Modeling}

\maketitle

\section*{Abstract}

%\begin{abstract}

Poroelastic wave simulations are important for many applications relating fluid flow and wave characteristics in porous rock formations.
Reciprocity is a key physical property of wave propagation in porous media that is important for such applications, even when viscous dissipation is present. However, numerical poroelastic simulations often fail to reproduce reciprocal responses because the discretization does not preserve the balance between reversible wave dynamics and irreversible fluid–solid drag. To address this, we formulate the Biot equations in terms of a continuous evolution operator split into a reversible (skew-adjoint) wave part and an irreversible (self-adjoint, non-positive) Darcy part, including the leading-order Johnson--Koplik--Dashen correction. This structure clarifies why reciprocity holds in the continuous equations and how it is easily broken in discrete form.

Guided by this interpretation, we construct a symmetric second-order Strang-splitting scheme with half-step source injection. The method conserves energy in the reversible subsystem, treats Darcy dissipation unconditionally stably, and retains Courant limits similar to elastic solvers. Using a staggered pseudo-spectral discretization, we model multimode propagation in 2D VTI media and obtain cross-component reciprocity with a relative $L_2$ misfit below $2.5\times10^{-6}$, demonstrating that the discrete scheme inherits the symmetry properties of the continuous evolution operator.

\section{Introduction}

Wave propagation in fluid-saturated porous rocks is essential for seismic reservoir characterization, full-waveform inversion, and microseismic analysis. Biot's theory predicts fast-P, slow-P, and S-wave modes in anisotropic porous media (Biot, 1956a,b; Carcione, 1996). Reviews on computational poroelasticity summarize classical numerical approaches (Carcione and Goode, 1995; Carcione et al., 2010, Carcione, 2022). The most general poroelastic models typically include the full JKD formulation (Johnson et al., 1987) in both time domain (Stokke, 2024) and frequency domain (Jakobsen et al., 2025). In this work, we use the anisotropic Biot equations in evolutionary form with the lowest-order JKD correction and no memory variables, which adequately captures seismic-frequency attenuation and dispersion (Carcione, 1996; Masson et al., 2006). 

This paper examines source--receiver reciprocity in computational poroelasticity and accommodates moment-tensor, vectorial, and scalar pressure sources.  Moment-tensor sources, that are important for induced seismicity and microseismic modelling, remain understudied in poroelastic simulations (Tohti et al., 2021; Zheng et al., 2024), making their inclusion a novel aspect of this work. 
Reciprocity follows from linearity, causality, time invariance, and Onsager symmetry of the constitutive operators, and manifests as symmetry of Green's functions and interchangeability of source and receiver positions (de~Hoop, 1995; Wapenaar and Fokkema, 2004).
 Since Biot's equations share this structure, reciprocity carries over to poroelasticity, although only a limited number of studies address it explicitly (Sidler et~al., 2013), and even fewer if any consider moment-tensor sources.

Operator splitting provides a concise framework for analysing these symmetry properties. Classical splitting or partition methods separates reversible elastic--inertial dynamics from irreversible Darcy dissipation (Carcione, Morency, and Santos, 2010). In a semigroup setting, first-order Lie splitting (Lie, 1970) and second-order symmetric Strang splitting (Strang, 1968) render the reversible operator skew-adjoint and the dissipative operator self-adjoint and contractive (Pazy, 1983). This adjoint structure is preserved under symmetric time injection of both vectorial and moment-tensor sources.

Because reciprocity must also hold at the discrete level (Sidler et al, 2013), we employ a staggered-grid pseudospectral scheme (Carcione et al., 2010) that maintains the adjoint pairing of divergence and gradient operators. Combined with symmetric Strang splitting and time-symmetric application of vector and moment-tensor sources, the resulting propagator preserves adjoint symmetry - and thus reciprocity -  to $\mathcal{O}(\Delta t^2)$ (Sidler et al., 2013). Simulations in layered VTI media demonstrate stable multimode propagation, and cross-component reciprocity tests show errors nearly approaching machine precision.

\section{Time evolution of the dynamic Biot system}

%We derive a first-order time-evolution system for wave propagation in %general
%poroelastic media. Biot's seminal works (1956a,b; 1962) introduced a
%displacement-based formulation involving the relative solid–fluid %displacement
%$\mathbf{w}$, forming the basis of modern poroelastic wave theory. Later
%developments, including the velocity–stress formulations reviewed by %Carcione
%and collaborators, recast the equations into a form well suited for %explicit
%time integration. Following this framework, we combine the constitutive
%relations, momentum balance, and dynamic Darcy law into a closed first-%order
%system for stresses, pore pressure, solid velocity, and fluid flux, providing a
%natural starting point for operator partitioning and Lie–Strang splitting
%methods.

\subsection{Constitutive relations}

We adopt Biot's poroelasticity in the low-frequency limit (Carcione, 1996). The infinitesimal strain-rate tensor associated with the solid $\mathbf{v}$ is 
$
\boldsymbol{\varepsilon}(\mathbf{v}) = \tfrac12\bigl(\nabla\mathbf{v} + \nabla\mathbf{v}^{\top}\bigr)
$.
Here $\mathbf{\sigma }$ and $p$ is the solid stress and pore fluid pressure, respectively. In rate form, the anisotropic stress--strain relation is
\begin{equation}
\partial_t \boldsymbol{\sigma}
=
\mathbf{C} : \boldsymbol{\varepsilon}(\mathbf{v})
-
\boldsymbol{\alpha}\,\partial_t p,
\label{eq:stressrate}
\end{equation}
where $\mathbf{C}$ is the drained TI stiffness tensor and 
$\boldsymbol{\alpha}=\mathrm{diag}(\alpha_1,\alpha_1,\alpha_3)$ is the Biot effective-stress tensor.

Fluid mass conservation couples volumetric strain, pore pressure, and Darcy flux,
\begin{equation}
\frac{1}{M}\,\partial_t p
+
\boldsymbol{\alpha} : \boldsymbol{\varepsilon}(\mathbf{v})
+
\nabla\!\cdot \mathbf{q}
=
s,
\label{eq:mass}
\end{equation}
with $M$ the Biot modulus, $\mathbf{q}$ the fluid flux relative to the solid matrix, and $s$ a possible volumetric source. Equations \eqref{eq:stressrate}--\eqref{eq:mass} provide the quasistatic constitutive closure required for the momentum balance and dynamic Darcy law in the subsequent section.

\clearpage

\subsection{Momentum balance and dynamic Darcy law}

Dynamic coupling between the solid and fluid phases follows mixture momentum conservation,
\begin{equation}
\nabla\!\cdot\boldsymbol{\sigma} + {\bf f}
=
\rho\,\partial_t\mathbf{v}
+
\rho_f\,\partial_t\mathbf{q},
\label{eq:momentum}
\end{equation}
where $\rho $ is the bulk density, $\rho _{f}$ is the fluid density, $\mathbf{q}$ is the Darcy flux, and $\mathbf{f}$ is a body force source term. The above equation is a compact tensor form of the componentwise relations (18)--(19) in Carcione (1996).
In the low-frequency regime, the anisotropic dynamic Darcy law is
\begin{equation}
-\nabla p
=
\rho_f\,\partial_t\mathbf{v}
+
\mathbf{m}\cdot \,\partial_t\mathbf{q}
+
\mathbf{D}\cdot \,\mathbf{q}.
\label{eq:darcy}
\end{equation}
Here $\mathbf{m}=\mathrm{diag}(m_1,m_3)$ with $m_i = T_i\rho_f/\phi$ 
 where $\phi $ is the porosity, $T_{i}$ is the tortuosity in the $i$th direction and $\phi $ is the porosity. The inverse fluid mobility tensor $\mathbf{D} = \eta\,\mathbf{K}_0^{-1}$, where $\mathbf{K}$ is the permeability tensor and $\eta $ is the fluid viscosity. This expression matches the componentwise forms in equations (25)--(26) of Carcione (1996). 

In equation (4), the inertial coupling term 
$\mathbf{m}\,\partial_t\mathbf{q}$ is the leading-order dynamic correction in the low-frequency limit of the Johnson--Koplik--Dashen (JKD) permeability model (Johnson et al., 1987; Masson et al., 2006). Although small at seismic frequencies, it is retained because it is intrinsic to Biot's dynamic Darcy law and is implicitly included in the Kelvin--Voigt-type viscous approximation commonly used in numerical formulations (Carcione, 1996).

In implicit schemes such as the Crank--Nicolson discretization of Carcione and Goode (1995) and the operator-splitting strategy of Carcione (1996), unconditional stability stems from the A-stable, semi-implicit time integrator rather than from any stabilizing property of the inertial term. In contrast, Masson et al.\ (2006) show that for a fully explicit staggered-grid velocity--stress scheme the inertial Darcy term plays a decisive numerical role: omitting $\mathbf{m}\,\partial_t\mathbf{q}$ produces unconditional instability, whereas retaining it restores a standard Courant-type CFL condition analogous to elastic wave propagation. Thus, even if physically small, the inertial term is essential both for preserving the correct low-frequency limit of dynamic Darcy flow and for ensuring stability in explicit poroelastic time-stepping schemes.

\subsection{First-order evolution equations}

Starting from the constitutive relations, mass balance, momentum balance,
and the dynamic Darcy law introduced in the preceding subsections, we
collect the evolution equations for the state variables
$(\boldsymbol{\sigma},\,p,\,\mathbf{v},\,\mathbf{q})$ in explicit
first-order form suitable for structure-preserving time integration.

The stress and pressure updates follow directly from the constitutive
relations and fluid mass conservation:
\begin{equation}
\partial_t \boldsymbol{\sigma}
=
\mathbf{C} : \boldsymbol{\varepsilon}(\mathbf{v})
-
\boldsymbol{\alpha}\,\partial_t p,
\label{eq:ev_sigma_final}
\end{equation}
\begin{equation}
\partial_t p
=
M\left(
s
-
\boldsymbol{\alpha} : \boldsymbol{\varepsilon}(\mathbf{v})
-
\nabla\!\cdot \mathbf{q}
\right).
\label{eq:ev_p_final}
\end{equation}

To obtain evolution equations for $\mathbf v$ and $\mathbf q$,
we invert the coupled inertial block formed by the momentum balance and
dynamic Darcy law.  Introducing the diagonal inertial operator
\[
\boldsymbol{\Delta} = \rho\,\mathbf m - \rho_f^{\,2}\mathbf I ,
\]
the velocity and flux updates become
\begin{equation}
\partial_t \mathbf v
=
\boldsymbol{\Delta}^{-1}\!\cdot
\left[
\mathbf m\cdot(\nabla\!\cdot\boldsymbol{\sigma} + \mathbf f)
+ \rho_f\,\nabla p
+ \rho_f\,\mathbf D\cdot\mathbf q
\right],
\label{eq:ev_v_final}
\end{equation}
\begin{equation}
\partial_t \mathbf q
=
\boldsymbol{\Delta}^{-1}\!\cdot
\left[
-\rho_f\,(\nabla\!\cdot\boldsymbol{\sigma} + \mathbf f)
- \rho\,\nabla p
- \rho\,\mathbf D\cdot\mathbf q
\right].
\label{eq:ev_q_final}
\end{equation}

Equations \eqref{eq:ev_sigma_final}--\eqref{eq:ev_q_final} constitute the
complete homogeneous first-order Biot system advanced by the
structure-preserving splitting scheme.  The source terms $s$ and
$\mathbf f$ enter only through the physically appropriate channels
(mass balance and momentum balance), ensuring consistency with the
continuous Biot model.

\section{Structure-preserving Time Integration}

The evolution equations (5--8) form a first-order system for the solid
velocity $\mathbf v$, fluid flux $\mathbf q$, pore pressure $p$, and
stress tensor $\boldsymbol{\sigma}$. If we define the state vector 
\[
\Psi = (\mathbf v,\,\mathbf q,\,p,\,\boldsymbol{\sigma})^\top,
\]
we can express the equation equations (5-8) in block matrix form:
\[
\partial_t \Psi = H\Psi + F(t).
\]
The above equation reminds us about the Schroedinger equation in quantum mechanics. 
Following the evolution operator approach in quantum mechanics, we express the solution to the above equation formally as 
\[
\Psi(t+\Delta t) = e^{\Delta t H}\,\Psi(t).
\]
The operator $H$ (that reminds us about the Hamiltonian operator in quantum mechanics) 
contains the full coupled poroelastic physics. At the
continuum level, it combines two qualitatively different mechanisms:
poroelastic wave propagation, arising from elastic and inertial
couplings in the solid–fluid mixture, and viscous Darcy relaxation, due
to drag on the relative fluid flux. These mechanisms have fundamentally
different physical effects on the poroelastic energy. Wave propagation
is reversible and conserves energy, whereas Darcy flow is irreversible
and dissipative. It is therefore natural—both physically and
mathematically—to separate these contributions so that they can be
treated according to their distinct properties.

This motivates decomposing the Biot operator as
\[
H = T + V,
\]
where $T$ collects the elastic, inertial, and stress–strain couplings
responsible for reversible wave motion, and $V$ contains the symmetric
Darcy drag acting on the relative flux. This decomposition is analogeous to the decomposition of the Hamiltonian operator in quantum mechanics into kinetic and potential energy operators. In our poroelastic context, this split is meaningful
because the two parts $T$ and $V$ have complementary mathematical properties in the
poroelastic energy inner product: $T$ is skew-adjoint and generates a
unitary, non-dissipative evolution (Pazy, 1983; Kato, 1995), while $V$
is self-adjoint and non-positive and generates a contractive, monotone
decay of energy. In the
poroelasticity literature these are often described as the non-stiff and
stiff components (Carcione et al., 2010), but the terminology
“reversible’’ and “irreversible’’ more precisely reflects their
mathematical structure and clarifies why stiffness arises exclusively
from the dissipative sector.

Substituting $H = T + V$ into the exact evolution operator gives
\[
\Psi(t+\Delta t)
=
e^{\Delta t(T+V)}\,\Psi(t),
\]
but the operators $T$ and $V$ do not commute,
\[
[T,V]\neq 0,
\]
and so the exponential $e^{\Delta t(T+V)}$ cannot be evaluated in closed
form. 

Operator-splitting methods approximate this exponential by
compositions of the subflows generated individually by $T$ and $V$
(Trotter, 1959; Strang, 1968). The Lie product formula,
\[
S_{\Delta t}^{(1)}
=
e^{\Delta t T}\,e^{\Delta t V},
\]
provides a first-order approximation, while the symmetric Strang formula
\[
S_{\Delta t}^{(2)}
=
e^{\frac{\Delta t}{2}V}\,
e^{\Delta t T}\,
e^{\frac{\Delta t}{2}V}
\]
achieves second-order accuracy. These approximations inherit structure
from the exact subflows because they are constructed directly from the
natural generators of the reversible and irreversible dynamics.

To understand why, we examine the properties of the subflows. Since $T$
is skew-adjoint, it generates a unitary group
\[
\|e^{tT}\Psi\|_E = \|\Psi\|_E,
\]
corresponding to perfectly reversible poroelastic wave propagation.
Since $V$ is self-adjoint and non-positive, it generates a contractive
$C_0$ semigroup
\[
\|e^{tV}\Psi\|_E \le \|\Psi\|_E,
\qquad t\ge0,
\]
representing irreversible viscous relaxation. Because Strang splitting
uses these exact subflows, it preserves the continuum operator structure
to $\mathcal{O}(\Delta t^2)$: energy is conserved in the $T$-step, Darcy
dissipation is monotone in the $V$-step, and no artificial damping of
fast-$P$, slow-$P$, or shear modes is introduced (Carcione, 1996).

Furthermore, the Strang propagator satisfies the adjoint symmetry
\[
\left(S_{\Delta t}^{(2)}\right)^*
=
S_{-\Delta t}^{(2)}
+ \mathcal{O}(\Delta t^2),
\]
ensuring reciprocity and the correct adjointness relations of the
discrete evolution (Sidler et al., 2013). These
properties arise directly from the semigroup structure of the
generators, without requiring additional stabilization such as the
inertial JKD correction used in monolithic schemes (Masson et al., 2006).

Source terms lie outside both generators and are therefore added separately
using a symmetric half-kick: each force or moment is applied in two equal
impulses, one before and one after the conservative update. This Strang‑consistent,
time‑symmetric placement preserves the adjoint structure needed for reciprocity
while leaving the $T$‑ and $V$‑subflows unchanged.

\section{Spatial Discretization}

The structure-preserving Lie and Strang integrators presented here 
are independent of the method used for spatial
discretization, but for the numerical examples we use a staggered-grid
pseudo-spectral method. Solid velocities and fluid fluxes are stored at
staggered locations relative to stresses and pressures, and spatial derivatives
are evaluated spectrally using FFTs. In Fourier space, a derivative is computed as (Carcione et al., 2010)
\[
\partial_x u \;\approx\; \mathcal{F}^{-1}\!\left( ik_x\,\widehat{u} \right),
\]
with appropriate phase factors applied to enforce the staggering. Because the
multipliers $ik_x$ and $ik_z$ are purely imaginary, the discrete gradient and
divergence remain adjoint in the discrete energy inner product, ensuring that the
semi-discrete $T$- and $V$-blocks retain the same skew-adjoint and
self-adjoint negative-semidefinite structure as in the continuous formulation.

All examples use a uniform $N_x \times N_z$ grid with spacings $\Delta x$ and
$\Delta z$, and a time step $\Delta t$ chosen from the Courant condition of the
conservative block. A cosine-taper absorbing layer surrounds the computational
domain, and a two-thirds dealiasing rule (or a mild exponential filter)
suppresses spectral wrap-around. Material parameters and source characteristics
are given in the figure caption.

\section{Numerical Examples}

We consider a 2D layered VTI poroelastic model with identical material
parameters in all experiments (Figure 1).  The only changes are the source type and 
the saturating fluid.  This allows us to isolate the effects of viscosity, 
permeability, and shear‐wave excitation while keeping the geometry fixed.

%\paragraph{Example 1: Monopole source, oil saturation.}
We first use an isotropic monopole source that excites only P--waves.
The medium is saturated with a realistic reservoir oil.  In this case the 
slow P--wave is heavily attenuated, and the resulting wavefield is dominated 
by the fast P--mode and its reflections within the VTI layered structure.
This experiment serves as a reference for poroelastic wave propagation under
realistic viscous conditions.

%\paragraph{Example 2: Monopole source, gas saturation.}
Next, we use the \emph{same} monopole source and the \emph{same} VTI model, 
but replace the oil with a low‐viscosity gas.  
The reduced viscous drag and increased hydraulic diffusivity allow the slow 
P--wave to become propagative.  In this experiment, both the fast and slow 
compressional modes are clearly observed, providing a direct comparison with 
Example~1 and illustrating the strong dependence of Biot's slow mode on the 
fluid properties.

%\paragraph{Example 3: Moment‐tensor source, oil saturation.}
Finally, we return to the oil‐saturated medium of Example~1 but replace the 
monopole excitation with a moment‐tensor source containing off‐diagonal 
components.  This generates both P-- and S--waves, enabling us to examine 
wave conversions and anisotropic shear propagation in the same VTI model.
This experiment demonstrates that the numerical scheme can handle the full 
poroelastic wave system, including shear waves, mode conversions, and 
anisotropic coupling.

\plot{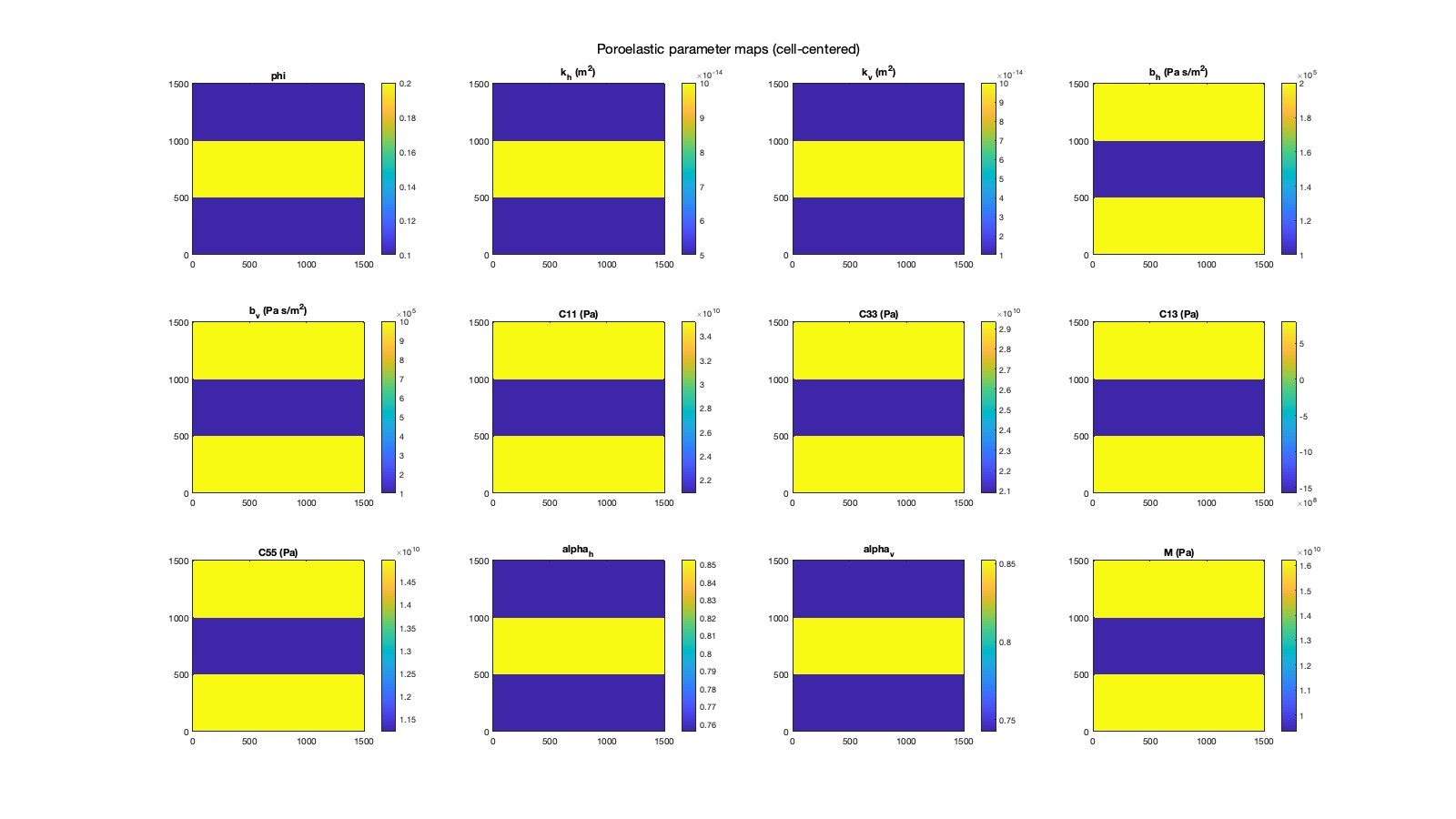}{width=\columnwidth}{
Heterogeneous VTI model.}

\plot{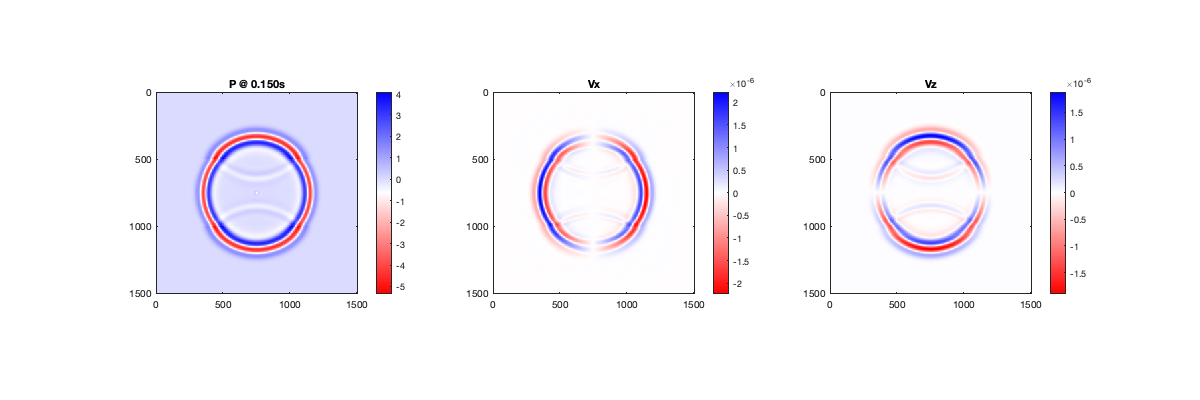}{width=\columnwidth}{
Snapshot of the wavefield in the layered VTI model shown in figure 1 for a monopole source ($M_{xx} = 1$, $M_{zz}=1$ and $M_{xz} = 0$) and oil saturation ($\eta = 1e-3$ Pa$\cdot s$).}

\plot{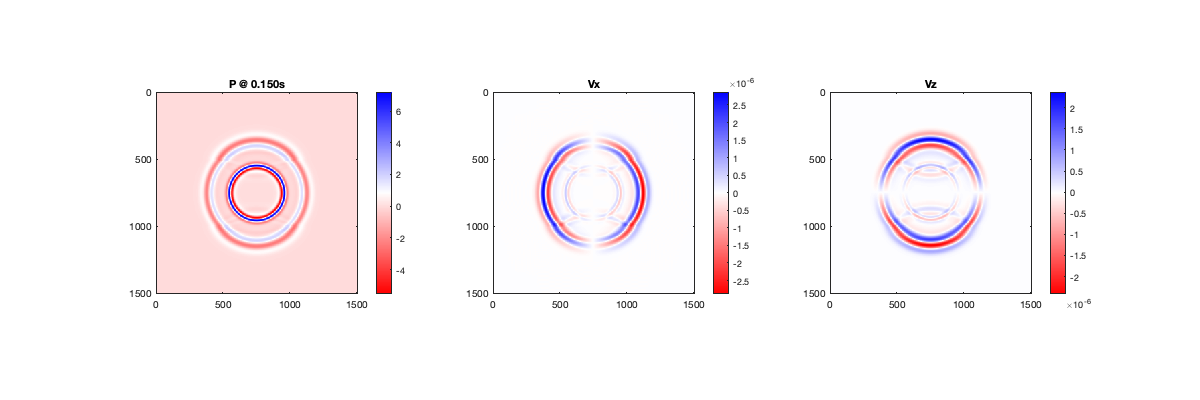}{width=\columnwidth}{
The same as Figure 2, but for gas saturation ($\eta = 1e-7$ Pa$\cdot s$).}

\plot{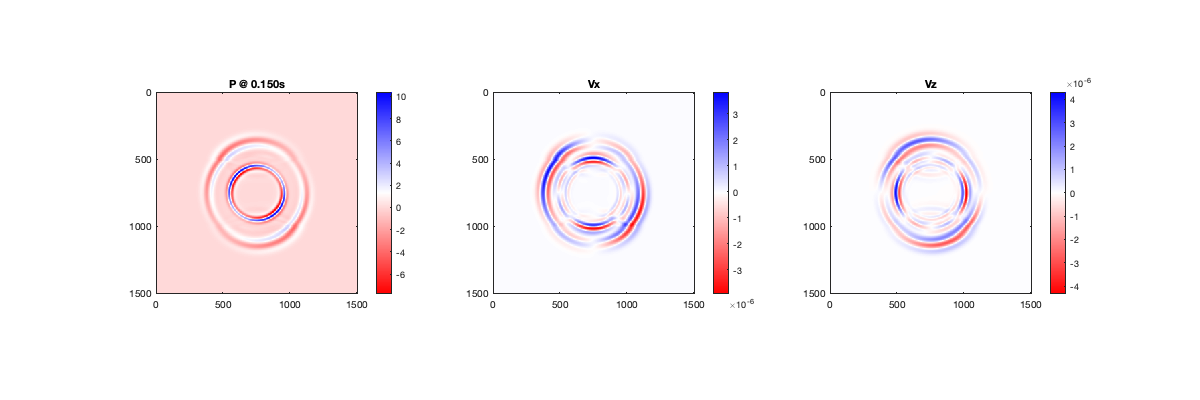}{width=\columnwidth}{
The same as in Figure 3 but for a moment tensor source the also generates S-waves 
($M_{xx} = 1$, $M_{zz}=1$ and $M_{xz} = 0.5$).}

\section{Discussion: Reciprocity}

A stringent verification of a poroelastic wave solver is to check whether it
preserves reciprocity, which follows from the mixed skew‑adjoint and
self‑adjoint structure of the continuous Biot equations. A discretization that
respects this structure should yield reciprocal signals, each given by the
convolution of a Green’s tensor component with the specified source‑time
function (Karpfinger et al., 2009). Even small violations of staggering, adjointness, or source injection
typically lead to observable non‑reciprocal behaviour.
We follow the cross‑component test of Sidler et al.\ (2013) using two
locations $\mathbf{r}_a$ and $\mathbf{r}_b$. A horizontal force at
$\mathbf{r}_a$ with vertical velocity recorded at $\mathbf{r}_b$ therefore gives
the signal corresponding to $G_{zx}(\mathbf{r}_b,\mathbf{r}_a,t)$ convolved with
the source‑time function, while interchanging source and receiver and swapping
force and measurement directions yields the corresponding $G_{xz}$ signal. For a
reciprocal medium these traces should coincide. Figure~5 shows an almost perfect
overlap for a viscosity of $10^{-3}$~Pa$\cdot$s, with a relative $L^2$ misfit
below $2.5\times10^{-6}$ in the reflection‑free window, demonstrating that the
discretization preserves the expected adjoint symmetry of the Biot system.

\renewcommand{\figdir}{Fig}

\begin{figure}[ht]
\centering
\includegraphics[height=4cm, width = 0.8\columnwidth]{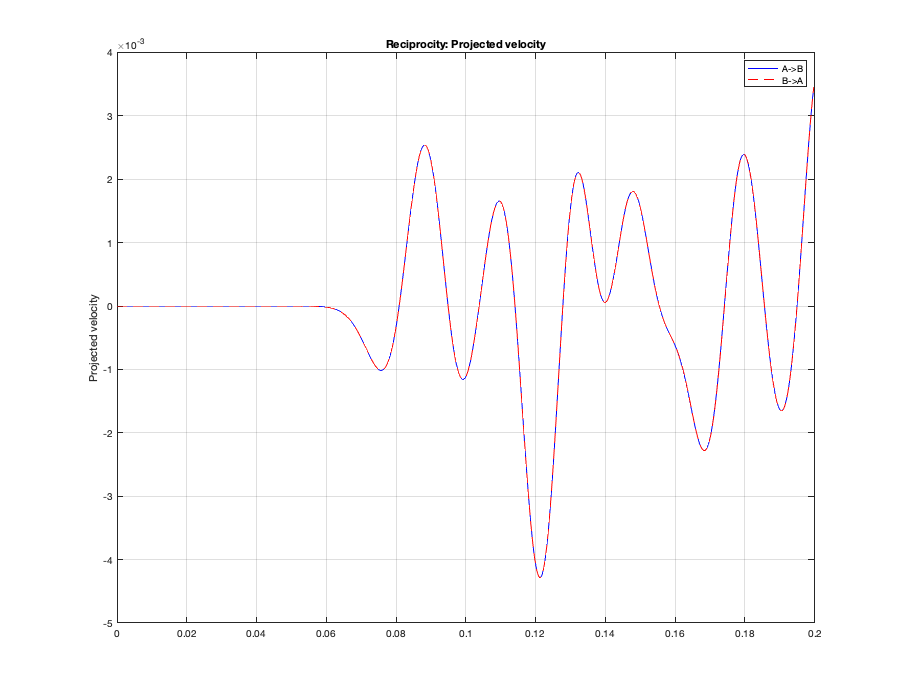}
\caption{Cross‑component reciprocity test. 
Blue: $G_{z x}(\mathbf{r}_b,\mathbf{r}_a,t)$. 
Red dashed: $G_{x z}(\mathbf{r}_a,\mathbf{r}_b,t)$.}
\label{fig:reciprocity}
\end{figure}

\section*{Concluding remarks}

A time-symmetric partitioning (Strang-splitting) update combined with symmetric
half-kick source injection preserves the adjoint symmetry required for
reciprocal Green’s functions in poroelastic wave simulations. Tests in layered
VTI media confirm stable multimode propagation in high-viscosity settings, with
reciprocity errors near machine precision.

For seismic frequencies, we employ the Biot--JKD model with only the
lowest-order dynamic-permeability correction, which captures attenuation,
dispersion, and stabilizing Darcy effects without the cost of full
fractional-derivative formulations.

At higher frequencies—such as in borehole acoustics and ultrasonic laboratory
applications—the full JKD model becomes necessary. Because the complete
formulation remains linear and causal when expressed via memory variables,
reciprocity should carry over if these variables are updated in an
adjoint-consistent manner, making this a natural direction for future work.

%\end{document}

\clearpage 
\section*{REFERENCES}

Biot, M. A., 1956a, Theory of propagation of elastic waves in a fluid-saturated porous
solid. I. Low-frequency range: \textit{Journal of the Acoustical Society of America}, 28, 168--178.

Biot, M. A., 1956b, Theory of propagation of elastic waves in a fluid-saturated porous
solid. II. Higher frequency range: \textit{Journal of the Acoustical Society of America}, 28, 179--191.

Carcione, J. M., and Quiroga-Goode, G., 1995, Some aspects of the physics and numerical modeling of Biot compressional waves, J. Comput. Acous., 3, 261-280. 

Carcione, J. M., 1996,
Wave propagation in anisotropic, saturated porous media:
Plane-wave theory and numerical simulation: \textit{J.\ Acoust.\ Soc.\ Am.}, 99, 2655--2666.

Carcione, J. M., C. Morency, and J. E. Santos, 2010, Computational poroelasticity: A review: \textit{Geophysics}, 75, 5, 75A229–75A243.

Carcione, J. M., 2022, \textit{Wave Fields in Real Media}, 4th ed.: Elsevier.

de Hoop, A. T., 1995, Handbook of Radiation and Scattering of Waves: Acoustic waves in fluids, elastic waves in solids, electromagnetic waves: Academic Press.

Jakobsen, M., J. S. Stokke, K. Kumar, and A. F. Radu, 2025,
Frequency-domain Biot--Allard equations for isotropic and anisotropic
poroelastic media: Two-field formulations and iterative splitting:
\textit{Research Square}, preprint,
https://doi.org/10.21203/rs.3.rs-8293824/v1.

Johnson, D. L., J. Koplik, and L. Dashen, 1987, Theory of dynamic permeability and
tortuosity in fluid-saturated porous media: \textit{Journal of Fluid Mechanics}, 176, 379--402.

Karpfinger, F., T. M. Müller, and B. Gurevich, 2009, Green’s functions and radiation patterns in poroelastic solids revisited: Geophysical Journal International, 178, 327–337.

Kato, T., 1995, \textit{Perturbation theory for linear operators}, Springer-Verlag.

Lie, S., 1970, The Lie product formula for semigroups generated by unbounded operators:
\textit{Mathematische Annalen}, 188, 235--242.

Masson, Y. J., S. R. Pride, and K. T. Nihei, 2006, Finite difference modeling of Biot’s
poroelastic equations at seismic frequencies: \textit{Journal of Geophysical Research}, 111, B10305.

Pazy, A., 1983, \textit{Semigroups of linear operators and applications to partial differential
equations}, Springer-Verlag.

Sidler, R., J. M. Carcione, and K. Holliger, 2013,
A pseudo-spectral method for the simulation of poro-elastic seismic wave propagation
in 2D polar coordinates using domain decomposition:
\textit{Journal of Computational Physics}, 235, 846--864.

Strang, G., 1968, On the construction and comparison of difference schemes:
\textit{SIAM Journal on Numerical Analysis}, 5, 506--517.

Stokke, J. S., M. Jakobsen, K. Kumar, and F. A. Radu, 2025, 
A history-dependent dynamic Biot model: In Sequeira et al. (Eds.), 
ENUMATH 2023, Lecture Notes in Computational Science and Engineering, 
154, 360–368, Springer.

Tohti, M., Y. Wang, W. Xiao, Q. Di, K. Zhou, J. Wang, S. An, and S. Liao, 2021,
Numerical simulation of seismic waves in 3-D orthorhombic poroelastic medium with
microseismic source implementation: \textit{Geophysical Journal International}, 227, 1012--1027.

Trotter, H. F., 1959, On the product of semi-groups of operators:
\textit{Proceedings of the American Mathematical Society}, 10, 545--551.

Wapenaar, K., and J. Fokkema, 2004, Reciprocity theorems for diffusion, flow, and waves: Journal of Applied Mechanics, 71(1), 145–150.

Zheng, J., T. Li, J. Xie, and Y. Sun, 2024, Improved numerical solution of anisotropic poroelastic wave equation in microseismicity: Graphic process unit acceleration and moment tensor implementation: \textit{Geophysical Prospecting}, 72(6)

\end{document}